\documentclass[12pt]{article}
\title{Mysterious duality and helical line bundles\\ on del Pezzo surfaces}
\author{Alastair King}
\date{11 July 2025} 

\usepackage[a4paper]{geometry}
\usepackage{amssymb,amsthm,array}
\usepackage[reqno]{amsmath}
\usepackage[usenames,dvipsnames]{xcolor}
\usepackage[colorlinks, linkcolor=RoyalBlue,anchorcolor=Periwinkle,
    citecolor=Orange,urlcolor=Emerald]{hyperref}
\usepackage{tikz}
  \usetikzlibrary{calc}

\theoremstyle{definition}
\newtheorem{theorem}{Theorem}[section]
\newtheorem{proposition}[theorem]{Proposition}

\theoremstyle{definition}

\newtheorem{example}[theorem]{Example}

\newtheorem{remark}[theorem]{Remark}

\numberwithin{equation}{section}
\numberwithin{figure}{section}

\renewcommand{\epsilon}{\varepsilon}

\renewcommand{\leq}{\leqslant}
\renewcommand{\geq}{\geqslant}

\newcommand{\ZZ}{\mathbb{Z}}

\newcommand{\Cbar}{\overline{C}}

\newcommand{\isom}{\cong}
\newcommand{\half}{{\textstyle\frac12}}
\newcommand{\sub}{\subseteq}
\newcommand{\spn}[1]{\left\langle #1 \right\rangle}

\newcommand{\lra}{\longrightarrow}

\newcommand{\cO}{\mathcal{O}}
\newcommand{\PP}{\mathbb{P}}
\newcommand{\FF}{\mathbb{F}}
\newcommand{\Cl}{\operatorname{Clf}}
\newcommand{\Lie}{\operatorname{Lie}}
\newcommand{\g}{\mathfrak{g}}
\newcommand{\dual}{^*}
\newcommand{\Etil}{\widetilde{E}}

\renewcommand{\sl}[1]{\mathfrak{sl}_{#1}}
\newcommand{\gl}[1]{\mathfrak{gl}_{#1}}
\newcommand{\so}[1]{\mathfrak{so}_{#1}}
\newcommand{\su}[1]{\mathfrak{su}_{#1}}
\newcommand{\uu}[1]{\mathfrak{u}_{#1}}
\newcommand{\ee}[1]{\mathfrak{e}_{#1}}

\newcommand{\rtsys}{\Phi}
\newcommand{\gurep}[2]{R_{#1}(#2)}
\newcommand{\gurepdual}[2]{R_{#1}(#2)\dual}
\newcommand{\extpow}[2]{\Lambda^{#1}(#2)}

\begin{document} \maketitle

\begin{abstract}
Mysterious duality is a relationship, described by Vafa in 2000, 
between $\textstyle\frac12$-BPS branes in Type II supergravity in dimension $D=d+2$ 
and rational curves on del Pezzo surfaces of degree $d$.
We show that both sides of this correspondence can be linked to a $\mathbb{Z}_d$ grading of the 
Lie algebra~$E_8$.
In addition, we show that the relevant rational curves correspond to 
`helical' line bundles, that is, line bundles that can appear in a helix on the del Pezzo surface. 
\end{abstract}

\section{Mysterious duality}\label{sec:MD}

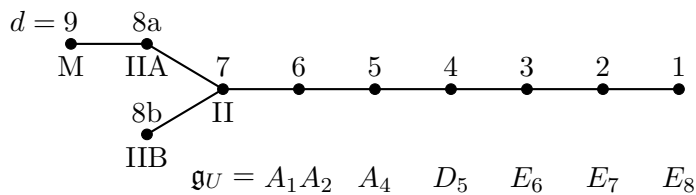
\begin{figure}[h]\centering
\begin{tikzpicture}
\newcommand{\down}{(0,-0.3)}
\newcommand{\upit}{(0,0.3)}
\foreach \j in {1,2,3,4,5,6,7} {\coordinate (a\j) at (9 - \j,0);}
\draw (a7)+\down node {\small II};
\coordinate (a8a) at (1,0.6); \draw (a8a)+\down node {\small IIA};
\coordinate (a8b) at (1,-0.6); \draw (a8b)+\down node {\small IIB};
\coordinate (a9) at (0,0.6); \draw (a9)+\down node {\small M};
\foreach \t/\h in {a1/a2, a2/a3, a3/a4, a4/a5, a5/a6, a6/a7, a7/a8a, a7/a8b, a8a/a9} 
  {\draw[thick] (\t)--(\h);} 
\foreach \j in {1,2,3,4,5,6,7,8a,8b,9} 
{ \filldraw[black] (a\j) circle(2pt);
 \draw (a\j)+\upit node {\small \j};}
 \draw (a9)+(-0.5,0.3) node {\small $d=$};
 \foreach \j/\D in {1/E_8, 2/E_7, 3/E_6, 4/D_5, 5/A_4, 6/A_1A_2} 
 \draw (9-\j,-1.2) node {\small $\D$};
 \draw(2,-1.2) node {\small $\g_U = $};
\end{tikzpicture}
\caption{del Pezzo/supergravity tree}
\label{fig:dPtree}
\end{figure}

It is an observation going back to Julia \cite{Jul81} that the classification of maximal, 
or Type~II, supergravity theories matches the classification of del Pezzo surfaces, 
via the appearance of an extended $E$ series of Lie groups (see Figure~\ref{fig:dPtree}).
In Type II supergravity in dimension $D$, the split real form group with Lie algebra~$\g_U$ 
acts on the scalar fields of the theory. 
In the corresponding toroidal compactification of string theory, the integral form of the group 
appears as the $U$-duality symmetry~\cite{HuTo95}.

For a del Pezzo surface of degree $d=D-2$, 
the intersection lattice naturally carries an action of the Weyl group of $\g_U$
(see e.g. \cite[\S23-26]{Ma}).
Note: $\g_U$ is defined for every~$d$, but not all are listed here 
(see Figure~\ref{fig:Etil-grading} later for the complete list).

There is one deformation type of degree $d$ del Pezzo surfaces $dP_d$ for each $1\leq d\leq 9$, 
except for $d=8$ where we have the two ruled surfaces $\FF_0=\PP^1\times\PP^1$ and $\FF_1=\PP^2(1)$,
i.e.~$\PP^2$ blown up at one point.
We shall denote these surfaces by $dP_{8a}=\FF_1$ and $dP_{8b}=\FF_0$ and recall that $dP_9=\PP^2$ .
The surfaces of lower degree are obtained by further blowing up either $dP_{8a}$ or $dP_{8b}$.

Matching this behaviour, there is one Type~II supergravity theory for each $3\leq D\leq 11$, 
except for $D=10$, where there are two, known as IIA and IIB, which are the low energy limits of the corresponding string theories.
One of these (IIA) arises from compactification on $S^1$ of the unique theory in $D=11$, 
which is expected to be the low energy limit of M theory. 
Both theories give rise to the same theories in $D\leq 9$ on further $S^1$-compactification.

The correspondence was christened `mysterious duality'
by C.~Vafa in 2000 (\cite{Vaf00, IqNeVa01}) based on the further observation that
the intersection lattice of $dP_d$ can be identified 
with the charge lattice for $\half$-BPS branes in Type II string theory in dimension $D=d+2$ 
in such a way that the classes of the branes match with the classes of certain rational curves.
Note that the addition of one to the rank of the intersection lattice on blowing up matches with the additional winding number for the branes around the new $S^1$ on compactification.

In its simplest incarnation, mysterious duality makes the connection between 
M-theory or $D=11$ supergravity, which contains 
electro-magnetic (E-M) dual $2$-branes and $5$-branes,
and the projective plane $\PP^2$, a surface of
degree $d=9$ in its anti-canonical embedding, 
containing rational curves $C_3$ (lines) and $C_6$ (conics)
of degrees 3 and 6, which are dual in the sense that
\begin{equation}
\label{eq:EMC}
C_3\cup C_6 \in |-K|,
\end{equation}
the anti-canonical linear system.

To make the parallel numerology a bit clearer, observe that the 
supergravity theory in $D=11$ contains 3-form and 6-form potentials $A_3$ and $A_6$
which are E-M dual in the sense that 
\begin{equation}
\label{eq:EMA}
dA_6=*dA_3+A_3\wedge dA_3
\end{equation}
The $M2$-brane worldsheet corresponds to a 
3-dimensional singularity/source in the supergravity, 
which carries a charge given by integrating $*dA_3$ over a 7-sphere 
linking the worldsheet. 
The fact that $D=d+2$ is crucial in matching the numerology ($3+6=9$)
in \eqref{eq:EMC} and the leading term of \eqref{eq:EMA}.

As we will explain in this article, both sides of this somewhat tenuous correspondence
can also be linked to the Lie algebra decomposition
\begin{equation}
\label{eq:EMM}
 \ee{8}=\sl{9}\oplus\Lambda^3\oplus \Lambda^6.
\end{equation}

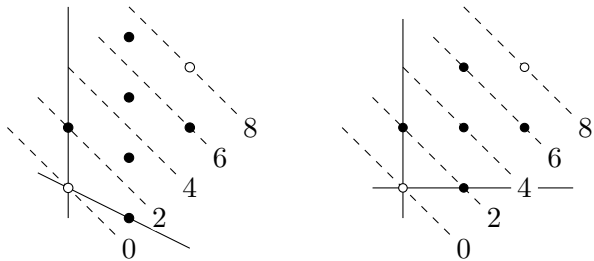
\begin{figure}[h]\centering
\begin{tikzpicture}[scale=0.8]
\draw (-0.5,0.25)--(2,-1) (0,-0.5)--(0,3);
\draw[dashed] (-1,1)--(1,-1) (-0.5,1.5)--(1.5,-0.5) (0,2)--(2,0) (0.5,2.5)--(2.5,0.5) (1,3)--(3,1);
\foreach \x/\y in {0/1, 1/-0.5, 1/0.5, 1/1.5, 1/2.5, 2/1}
{ \filldraw[black] (\x,\y) circle(2.2pt); }
\foreach \x/\y in {0/0, 2/2}
{ \draw[fill=white] (\x,\y) circle(2.2pt); }
\foreach \x/\y/\d in {1/-1/0, 1.5/-0.5/2, 2.0/0/4, 2.5/0.5/6, 3/1/8}
{ \filldraw[white] (\x,\y) circle(6pt);
  \draw (\x,\y) node {\small $\d$}; }
\end{tikzpicture}
\qquad
\begin{tikzpicture}[scale=0.8]
\draw (-0.5,0)--(2.8,0) (0,-0.5)--(0,2.8);
\draw[dashed] (-1,1)--(1,-1) (-0.5,1.5)--(1.5,-0.5) (0,2)--(2,0) (0.5,2.5)--(2.5,0.5) (1,3)--(3,1);
\foreach \x/\y in {1/0, 0/1, 1/1, 1/2, 2/1}
{ \filldraw[black] (\x,\y) circle(2pt); }
\foreach \x/\y in {0/0, 2/2}
{ \draw[fill=white] (\x,\y) circle(2pt); }
\foreach \x/\y/\d in {1/-1/0, 1.5/-0.5/2, 2.0/0/4, 2.5/0.5/6, 3/1/8}
{ \filldraw[white] (\x,\y) circle(6pt); ;
  \draw (\x,\y) node {\small $\d$}; }
\end{tikzpicture}
\caption{Intersection lattices of $dP_{8a}$ and $dP_{8b}$}
\label{fig:int-lat}
\end{figure}

The second (more convincing) example of mysterious duality 
involves type IIA and IIB supergravities in dimension $D=10$ 
in parallel with the two ruled del Pezzo surfaces $dP_{8a}=\FF_1$ and $dP_{8b}=\FF_0$.
Figure~\ref{fig:int-lat} shows their (rank 2) intersection lattices $I_X=H^2(X;\ZZ)$, marking 
(in white) the classes $0$ and $-K$ and (in black) the classes of rational curves $C$ 
which admit `dual' rational curves $\Cbar$, in the sense that $C+ \Cbar=-K$. 

The corresponding branes are known as F1 
(the fundamental string, corresponding the class of a fibre) 
and NS5 (in both IIA and IIB) 
and D0, D2, D4, D6 in IIA or D1, D3, D5 in IIB.

On the supergravity side, we begin with a small digression
and describe the (non-maximal) Type I theory in $D=10$
(cf.~\cite{Jul98}).
The analogue of \eqref{eq:EMM} in this case
is a more standard Lie algebra decomposition
\begin{equation} \label{eq:EMI}
 \so{16} = \gl{8} \oplus\Lambda^2\oplus \Lambda^6.
\end{equation}
This can be understood as a decomposition of either the complex form or the split real form,
which is often written $\so{8,8}$.
In particular, \eqref{eq:EMI} can help identify the Cartan involution $\sigma$ of the split real form,
whose maximal compact subalgebra (i.e.~the $+1$-eigenspace of $\sigma$) is $\so{8}\oplus \so{8}$.

On the right-hand side, we have the standard Cartan involution of $\gl{8}$ corresponding to a choice of metric,
giving the $\pm1$-eigenspace decomposition
\begin{equation}
\label{eq:gl8}
  \gl{8}=\so{8}\oplus [\Sigma^2_0\oplus 1], 
\end{equation} 
together with the Hodge-$*$ isomorphism $\Lambda^2\isom \Lambda^6$ coming from the metric,
giving an involution of $\Lambda^2\oplus \Lambda^6$.
Writing the eigenspaces of this involution as $\Lambda^{2|6}_{+}$ and~$\Lambda^{2|6}_{-}$,
the $\pm1$-eigenspace decomposition for $\sigma$ is then
\begin{equation}
\label{eq:EMIex}
 \so{16}= \bigl[ \so{8} \oplus \Lambda^{2|6}_{+} \bigr] 
   \oplus \bigl[ \Sigma^2_0 \oplus 1\oplus \Lambda^{2|6}_{-}\bigr]
\end{equation}
The massless field content of the Type I theory is given by the 
irreducible representations in the $-1$-eigenspace of $\sigma$:
the graviton corresponds to $\Sigma^2_0$,
the (scalar) dilaton to the trivial summand $1$,
and the B-field to $\Lambda^{2|6}_{-}$.
The $\Lambda^{2|6}_{+}$ is interpreted as a second $\so{8}$.

Applying the same method to \eqref{eq:EMM} also yields the Cartan involution 
for the split real form of $\ee{8}$, giving the $\pm1$-eigenspace decomposition
\begin{equation}
\label{eq:EMMex}
  \ee{8}= \bigl[ \so{9} \oplus \Lambda^{3|6}_{+} \bigr] 
   \oplus \bigl[ \Sigma^2_0 \oplus \Lambda^{3|6}_{-}\bigr].
\end{equation}
The massless fields in $D=11$ supergravity are just the graviton and C-field,
which is a 3-form.

To obtain the Type IIA and IIB theories for $D=10$ from Type I,
we should `add D-branes', which amounts to writing 
$\ee{8}=\so{16}\oplus S^{\pm}$ 
and identifying the spin representations $S^{\pm}$ with odd/even forms,
to extend \eqref{eq:EMI} to either 
\begin{equation}
\label{eq:EMIIA}
 \ee{8} = \bigl[ \gl{8}\oplus \Lambda^2\oplus \Lambda^6 \bigr]
 \oplus \bigl[ \Lambda^1\oplus\Lambda^3\oplus\Lambda^5\oplus\Lambda^7\bigr],
\end{equation}
for Type IIA, or
\begin{equation}
\label{eq:EMIIB}
 \ee{8} = \bigl[ \gl{8}\oplus \Lambda^2\oplus \Lambda^6 \bigr]
 \oplus \bigl[ \Lambda^0\oplus\Lambda^2\oplus\Lambda^4\oplus\Lambda^6
 \oplus\Lambda^8 \bigr],
\end{equation}
for Type IIB. 
Again the Cartan involution for the split real form comes from a metric and its Hodge-$*$ isomorphisms.

In \eqref{eq:EMIIA}, we will write the central (scalar) component of $\gl{8}$ as $\uu{1}$, 
but note that this rank 1 abelian Lie algebra is totally non-compact, 
i.e.~the Cartan involution is $-1$.
In \eqref{eq:EMIIB}, the three scalars $\uu{1}$, $\Lambda^0$, $\Lambda^8$ 
can be combined into the $S$-duality algebra $\sl{2}$ and so we end up with
\begin{equation}
\label{eq:EMIIAB}
\begin{aligned}
\text{Type IIA:}\quad \ee{8} & =  \bigl[ \sl{8}\oplus \uu{1} \bigr]
 \oplus \Lambda^1 \oplus \Lambda^2 \oplus \Lambda^3
 \oplus \Lambda^5 \oplus \Lambda^6 \oplus \Lambda^7\\
\text{Type IIB:}\quad \ee{8} & =  \bigl[ \sl{8}\oplus \sl{2} \bigr] 
 \oplus \Lambda^2 \cdot(2)\oplus \Lambda^4 \cdot(1)\oplus \Lambda^6 \cdot(2) 
\end{aligned}
\end{equation}
Note that the degrees of the forms (and their multiplicities)
correspond exactly to the degrees of the rational curves in $\FF_1$ and
$\FF_0$ marked in Figure~\ref{fig:int-lat}.

We will see in what follows (Theorem~\ref{thm:main} and Remark~\ref{rem:myst-dual})
that \eqref{eq:EMM} and \eqref{eq:EMIIAB} generalise 
to Lie algebra decompositions related to the similar `mysterious' correspondences 
between rational curves on del Pezzo surfaces of degree $d$ 
and branes in Type II supergravity in dimension $D=d+2$.

\newpage
\section{The shadow of $E_8$ in supergravity}\label{sec:E8}

\newcommand{\sprod}{s}

A maximal supergravity theory in dimension $D=d+2$ carries the symmetry
of the Poincar\'e superalgebra 
\[
  \bigl[ \so{d+1,1} \oplus P \bigr] \oplus Q^{32}
\]
where $P$ is the $D$-dimensional translation algebra
(so the first $\oplus$ is semi-direct).
Further, $Q^{32}$ is the odd part of the super algebra, 
which is a real (but usually not irreducible) spin representation 
of $\so{d+1,1}$ and satisfies $[Q,Q]\subseteq P$, 
i.e.~has an invariant symmetric product $\sprod\colon \Sigma^2 Q \to P$. 
We will see shortly why $Q$ should have dimension 32
and thus $D\leq 11$.

The `massless degrees of freedom' are determined by
fixing a null momentum $p\in P\dual$ (i.e.~$p^2=0$) 
to reduce the symmetry to the `helicity' algebra $\so{d}$,
that is, the semisimple part of the stabiliser of $p$, 
together with an action of the Clifford algebra $\Cl(Q^{16})$,
where $Q^{16}$ is the natural quotient of $Q^{32}$,
determined by the degenerate quadratic form $p\circ\sprod$, 
and is a real spin representation of $\so{d}$.
The `massless supermultiplet' is the 256-dimensional irreducible representation of $\Cl(Q^{16})$. 
This decomposes as $(128)_B \oplus (128)_F$ under the even part of $\Cl(Q^{16})$, 
where the summands are identified as bosonic ($B$) or fermionic ($F$) by their spin as representations of $\so{d}$.

For a supergravity theory, 
the highest spin representation of $\so{d}$ occurring
in the bosonic part $(128)_B$ should be the spin 2 graviton, 
i.e.~the trace-free symmetric square $\Sigma^2_0(d)$.
This requires that the reduced $Q$ should be spanned by 8 (complex) spin-$\frac{1}{2}$ creation operators 
and 8 annihilation operators, explaining why $\dim Q=16$.

In the maximum case $d=9$ ($D=11$) the spin representation $Q^{16}$ 
is irreducible and unique, but when $d=8$ there are two options:
 $Q^{16}=(8)_+ \oplus (8)_-$, which is the restriction from $\so{9}$ to $\so{8}$, 
or $Q^{16}=2(8)_{\pm}$,  which is not a restriction.
Restricting further to $\so{7}$, these both give
$Q^{16}=2(8)$, reproducing precisely the distinctive behaviour
of the IIA and IIB theories.

Thus the (massless) representation theoretic content of supergravity involves 
a real spin representation $\so{d}\to \so{16}$ together with an embedding
of representations $\Sigma^2_0(d)\to (128)_B$.
These can be understood as the even and odd parts of an embedding of 
$\sl{d}=\so{d}\oplus \Sigma^2_0(d)$ into $\ee{8}=\so{16}\oplus(128)_B$, 
where the parity is given by the eigenvalues of the Cartan involutions for the split real forms.
However, a explicit role for the actual Lie algebra structure of $\ee{8}$ in supergravity 
seems unclear.

The remaining (bosonic) field content of each Type II supergravity theory is given by the
remaining $\so{d}$ representations in the $(128)_B$, 
and it turns out that these are always given by $m$-form `gauge' fields ($1\leq m\leq d/2$), 
with multiplicities which are representations of an emergent `U-duality' 
symmetry\footnote{the case of $d/2$-forms is actually a little more subtle than this} $\g_U$,
together with scalar fields, which form the non-compact part of $\g_U$
(cf. \cite{CJLP98a}).

One way to understand this picture is via a decomposition of $\ee{8}=\Lie(E_8)$.
This decomposition can be described as the components of a $\ZZ_d$-grading of 
$\ee{8}$, as in the following theorem.
Note that, while the $\ee{8}$ emerges as the full U-duality symmetry $\g_U$ when $d=1$
(cf. \cite{MS}),
it is effectively present in all dimensions.

\goodbreak
\begin{theorem} \label{thm:main}
For each $d$ in Figure~\ref{fig:dPtree},
the Lie algebra $\ee{8}$ has a $\ZZ_d$-grading
\begin{equation}
\label{eq:E8decomp}
 \ee{8} = \bigl[ \sl{d} \oplus \g_U\bigr] \oplus
 \bigoplus_{m=1}^{d-1} \extpow{m}{d}\otimes \gurep{m}{d},
\end{equation}
where $\g_U$ is the `$U$-duality' Lie algebra 
and $\gurep{m}{d}$ 
are representations of $\g_U$, with 
\[
  \gurepdual{m}{d} \isom \gurep{d-m}{d}.
\]
Further, $\gurep{m}{d}$ is irreducible and miniscule (i.e.~its weights are a single Weyl group orbit),
unless $d=7$ and $m=1$ or $6$ (see Figure~\ref{fig:Rmd} for dimensions).
\end{theorem}

\begin{proof}
Following Kac (see \cite[\S10.5]{Helg}, \cite[\S8.6]{Kac} or \cite{Vin76}),
$\ZZ_d$-gradings of a simple Lie algebra $\g$ are determined by co-weights on the extended Dynkin diagram of $\g$,
which, in the simpler (multiplicity-free) case amounts to crossing some nodes on the diagram.
Thus the gradings that we wish consider here are given by the data listed in Figure~\ref{fig:Etil-grading}.
They are equivalent to certain $\ZZ$-gradings of the Kac-Moody algebra $\Lie(\Etil_8)$ .

The type of the deg 0 part is read off from the diagam with the cross(es) removed,
with an extra factor of $\uu{1}$ when there are two crosses.
Thus the Lie algebra $\g_U$ can be read off from the right-hand (red) part of the diagram, 
i.e.~ignoring the left-hand (blue) part which gives $\sl{d}$. 
In the multiplicity-free cases, the deg 1 part is read off from the nodes adjacent to the crossed nodes 
(see \cite{Vin76}), with some modification in the cases $d=8a,8b,9$.

\begin{figure}\centering
\newcommand{\xx}{0.14}
\newcommand{\cc}{0.11}
\newcommand{\drawcross}[1]{
 \draw[very thick] (#1)+(-\xx,-\xx) -- (#1); 
 \draw[very thick] (#1)+(\xx,\xx) -- (#1); 
 \draw[very thick] (#1)+(\xx,-\xx) -- (#1); 
 \draw[very thick] (#1)+(-\xx,\xx) -- (#1); 
}
\begin{tabular}{cccc}
 & $\g_U$ & diagram & co-weight
\\ \\
$d=1$ & $E_8$ &
\begin{tikzpicture}[scale=0.75]
\coordinate (a0) at (6,2);
\foreach \j in {1,2,3,4,5,6,7,8} {\coordinate (a\j) at (\j,1);}
\foreach \t/\h in {a0/a6, a1/a2, a2/a3, a3/a4, a4/a5, a5/a6, a6/a7, a7/a8} 
  {\draw[thick] (\t)--(\h);} 
\drawcross{a1};
\foreach \j in {2,3,4,5,6,7,8,0} 
 {\filldraw[red] (a\j) circle(\cc);}
\end{tikzpicture}
& $\omega_{1}$
\\ \\
$d=2$ & $E_7$ &
\begin{tikzpicture}[scale=0.75]
\coordinate (a0) at (6,2);
\foreach \j in {1,2,3,4,5,6,7,8} {\coordinate (a\j) at (\j,1);}
\foreach \t/\h in {a0/a6, a1/a2, a2/a3, a3/a4, a4/a5, a5/a6, a6/a7, a7/a8} 
  {\draw[thick] (\t)--(\h);} 
 \drawcross{a2};
\foreach \j in {1} 
 {\filldraw[blue] (a\j) circle(\cc);}
\foreach \j in {3,4,5,6,7,8,0} 
 {\filldraw[red] (a\j) circle(\cc);}
\end{tikzpicture} 
& $\omega_{2}$
\\ \\
$d=3$ & $E_6$ &
\begin{tikzpicture}[scale=0.75]
\coordinate (a0) at (6,2);
\foreach \j in {1,2,3,4,5,6,7,8} {\coordinate (a\j) at (\j,1);}
\foreach \t/\h in {a0/a6, a1/a2, a2/a3, a3/a4, a4/a5, a5/a6, a6/a7, a7/a8} 
  {\draw[thick] (\t)--(\h);} 
 \drawcross{a3};
\foreach \j in {1,2} 
 {\filldraw[blue] (a\j) circle(\cc);}
\foreach \j in {4,5,6,7,8,0} 
 {\filldraw[red] (a\j) circle(\cc);}
\end{tikzpicture}
& $\omega_{3}$
\\ \\
$d=4$ & $D_5$ &
\begin{tikzpicture}[scale=0.75]
\coordinate (a0) at (6,2);
\foreach \j in {1,2,3,4,5,6,7,8} {\coordinate (a\j) at (\j,1);}
\foreach \t/\h in {a0/a6, a1/a2, a2/a3, a3/a4, a4/a5, a5/a6, a6/a7, a7/a8} 
  {\draw[thick] (\t)--(\h);} 
 \drawcross{a4};
\foreach \j in {1,2,3} 
 {\filldraw[blue] (a\j) circle(\cc);}
\foreach \j in {5,6,7,8,0} 
 {\filldraw[red] (a\j) circle(\cc);}
\end{tikzpicture}
& $\omega_{4}$
\\ \\
$d=5$ & $A_4$ &
\begin{tikzpicture}[scale=0.75]
\coordinate (a0) at (6,2);
\foreach \j in {1,2,3,4,5,6,7,8} {\coordinate (a\j) at (\j,1);}
\foreach \t/\h in {a0/a6, a1/a2, a2/a3, a3/a4, a4/a5, a5/a6, a6/a7, a7/a8} 
  {\draw[thick] (\t)--(\h);} 
\drawcross{a5};
\foreach \j in {1,2,3,4} 
 {\filldraw[blue] (a\j) circle(\cc);}
\foreach \j in {6,7,8,0} 
 {\filldraw[red] (a\j) circle(\cc);}
\end{tikzpicture}
& $\omega_{5}$
\\ \\
$d=6$ & $A_1A_2$ &
\begin{tikzpicture}[scale=0.75]
\coordinate (a0) at (6,2);
\foreach \j in {1,2,3,4,5,6,7,8} {\coordinate (a\j) at (\j,1);}
\foreach \t/\h in {a0/a6, a1/a2, a2/a3, a3/a4, a4/a5, a5/a6, a6/a7, a7/a8} 
  {\draw[thick] (\t)--(\h);} 
 \drawcross{a6};
 \foreach \j in {1,2,3,4,5} 
 {\filldraw[blue] (a\j) circle(\cc);}
\foreach \j in {7,8,0} 
 {\filldraw[red] (a\j) circle(\cc);}
\end{tikzpicture}
& $\omega_{6}$
\\ \\
$d=7$ & $A_1u_1$&
\begin{tikzpicture}[scale=0.75]
\coordinate (a0) at (6,2);
\foreach \j in {1,2,3,4,5,6,7,8} {\coordinate (a\j) at (\j,1);}
\foreach \t/\h in {a0/a6, a1/a2, a2/a3, a3/a4, a4/a5, a5/a6, a6/a7, a7/a8} 
  {\draw[thick] (\t)--(\h);} 
 \drawcross{a7};
 \drawcross{a0};
 \foreach \j in {1,2,3,4,5,6} 
 {\filldraw[blue] (a\j) circle(\cc);}
\foreach \j in {8} 
 {\filldraw[red] (a\j) circle(\cc);}
\end{tikzpicture}
& $\omega_7= \omega_{3'}+\omega_{4'}$
\\ \\
$d=8a$ & $u_1$ &
\begin{tikzpicture}[scale=0.75]
\coordinate (a0) at (6,2);
\foreach \j in {1,2,3,4,5,6,7,8} {\coordinate (a\j) at (\j,1);}
\foreach \t/\h in {a0/a6, a1/a2, a2/a3, a3/a4, a4/a5, a5/a6, a6/a7, a7/a8} 
  {\draw[thick] (\t)--(\h);} 
 \drawcross{a8};
 \drawcross{a0};
\foreach \j in {1,2,3,4,5,6,7} 
 {\filldraw[blue] (a\j) circle(\cc);}
\end{tikzpicture}
& $\omega_{8a}=2\omega_{3'}+\omega_{2'}$
\\ \\
$d=8b$ & $A_1$ &
\begin{tikzpicture}[scale=0.75]
\coordinate (a0) at (6,2);
\foreach \j in {1,2,3,4,5,6,7,8} {\coordinate (a\j) at (\j,1);}
\foreach \t/\h in {a0/a6, a1/a2, a2/a3, a3/a4, a4/a5, a5/a6, a6/a7, a7/a8} 
  {\draw[thick] (\t)--(\h);} 
\drawcross{a7};
\foreach \j in {1,2,3,4,5,6,0} 
 {\filldraw[blue] (a\j) circle(\cc);}
\foreach \j in {8} 
 {\filldraw[red] (a\j) circle(\cc);}
\end{tikzpicture}
& $\omega_{8b}=2\omega_{4'}$
\\ \\
$d=9$ &  &
\begin{tikzpicture}[scale=0.75]
\coordinate (a0) at (6,2);
\foreach \j in {1,2,3,4,5,6,7,8} {\coordinate (a\j) at (\j,1);}
\foreach \t/\h in {a0/a6, a1/a2, a2/a3, a3/a4, a4/a5, a5/a6, a6/a7, a7/a8} 
  {\draw[thick] (\t)--(\h);} 
\drawcross{a0};
\foreach \j in {1,2,3,4,5,6,7,8} 
 {\filldraw[blue] (a\j) circle(\cc);}
\end{tikzpicture}
& $\omega_9=3\omega_{3'}$
\\ \\
&  &
\begin{tikzpicture}[scale=0.75]
\coordinate (a0) at (6,2);
\foreach \j in {1,2,3,4,5,6,7,8} {\coordinate (a\j) at (\j,1);}
\foreach \t/\h in {a0/a6, a1/a2, a2/a3, a3/a4, a4/a5, a5/a6, a6/a7, a7/a8} 
  {\draw[thick] (\t)--(\h);} 
\foreach \j in {0,1,2,3,4,5,6,7,8} 
 {\filldraw[black] (a\j) circle(\cc);}
 \foreach \j in {1,2,3,4,5,6} 
 {\draw (a\j) node[below=3pt] {\j};}
 \draw (a7) node[below=3pt] {$4'$};
 \draw (a8) node[below=3pt] {$2'$};
 \draw (a0) node[right=3pt] {$3'$};
\end{tikzpicture}
\end{tabular}
\caption{Some $\ZZ_d$-gradings of $E_8$ (or $\ZZ$-gradings of $\Etil_8$).}
\label{fig:Etil-grading}
\end{figure}

To analyse these gradings more closely, we follow an approach that
dates back to Coble \cite{Cob} (see also \cite{Ma}).
Consider the Minkowski lattice $\ZZ^{1,9}$
with orthogonal basis $h,e_1,\dots,e_9$, where $h^2=1$ and $e_i^2=-1$.
Using the additional element 
\begin{equation}
  \omega_0=3h-(e_1+\dots+ e_9),
\end{equation}
the root lattice  $\Lambda(\Etil_8)$ can be identified with $\omega_0^\perp\sub \ZZ^{1,9}$
and then the root system (strictly, the real roots) of $\Etil_8$ is given by
\[
 \rtsys(\Etil_8) = \{ \alpha\in \omega_0^\perp : \alpha^2=-2 \}.
\]
Since $\omega_0$ is a null vector, the root system $\rtsys(\Etil_8)$ is invariant under
$\alpha\mapsto \alpha+\omega_0$ and the image  in the quotient lattice
$\Lambda(E_8)=\omega_0^\perp/\langle\omega_0\rangle$ is the root system $\rtsys(E_8)$ of $E_8$.

\newcommand{\sgrad}{\lambda}

Any co-weight, that is, a linear map $\sgrad\colon\Lambda(\Etil_8)\to\ZZ$, 
gives a $\ZZ$-grading\footnote{that is, an additive $\ZZ$-labelling of the roots}
of $\rtsys(\Etil_8)$, which induces a $\ZZ_d$-grading of $\rtsys(E_8)$, 
for $d=\sgrad(\omega_0)$,
and thus a $\ZZ_d$-grading of the Lie algebra $\ee{8}=\Lie(E_8)$.
The maps $\sgrad$ that we will consider are given by taking the inner product (in $\ZZ^{1,9}$) with
\begin{equation} \label{eq:om_d}
  \omega_d= 3h - (e_1+\dots+e_n),
\end{equation} 
for $n=9-d$, including $n=1$ for $d=8a$, together with
\begin{equation}
\label{eq:om_8b}
  \omega_{8b}= 4h - 2(e_1+e_2). 
\end{equation}
Note that $\omega_d^2=d$
and these $\omega_d$ give precisely the co-weights listed in Figure~\ref{fig:Etil-grading},
relative to the basis of simple roots $e_i-e_{i+1}$, for $i=1,\ldots, 8$, together with $h-(e_1+e_2+e_3)$.

We can now show that these gradings give the claimed decompositions of $\ee{8}$.
In the main case \eqref{eq:om_d}, i.e.~$d\neq 8b$, 
we write the obvious orthogonal direct sum decomposition
\[
   \ZZ^{1,9}=\ZZ^{1,n} \oplus \ZZ^d,
\]
where $n=9-d$ and the second summand is negative definite standard Euclidean,
spanned by $e_{n+1},\dots,e_9$.
With respect to this decomposition, we have 
\[ 
  \omega_0=\omega_d -\Delta_d,
  \quad\text{where $\Delta_d=e_{n+1}+\dots+e_9$,} 
\]
and, for any real root $\alpha$ in $\Lambda(\Etil_8)=\omega_0^\perp$, 
we can write $\alpha=\beta+\gamma$,
that is, with $\beta\in\ZZ^{1,n}$ and $\gamma\in\ZZ^d$.
Then 
\[ 
  m=\omega_d\cdot\alpha=\omega_d\cdot\beta=\Delta_d\cdot\gamma.
\]
To obtain all $E_8$ roots, we can consider just $\Etil_8$ roots $\alpha$
with $0\leq m < d$.

If $m=0$, then $\Delta_d^\perp\subseteq \ZZ^d$ is $\Lambda(\sl{d})$ and 
$\omega_d^\perp\subseteq\ZZ^{1,n}$ is a known incarnation (see e.g.~\cite[\S23-26]{Ma}) 
of the root lattice $\Lambda(\g_U)$, for $\g_U$ as listed in Figure~\ref{fig:Etil-grading}. 
These are both negative definite lattices, without vectors $v$ with $v^2=-1$,
so if $\alpha^2=-2$, then 
either $\beta=0$ and $\gamma^2=-2$, so $\gamma$ is a root of $\sl{d}$,
or $\gamma=0$ and $\beta^2=-2$, so $\beta$ is a root of $\g_U$.

If $\alpha=\beta+\gamma$, as above, is a (real) $\Etil_8$ root with $1\leq m\leq d-1$, then we claim:

\begin{enumerate}
\item[(A)] $\gamma=-e_I=-\sum_{i\in I} e_i$, for $I$ any $m$-element subset of
$\{n+1,\dots,9\}$.
Thus, after projection onto the weight lattice $\ZZ^d/\spn{\Delta_d}$ of $\sl{d}$,
the $[\gamma]$ are the weights of $\extpow{m}{d}$.
In particular, this means that $\gamma^2=-m$ and so $\beta^2=m-2$.
\item[(B)] after projection onto $\ZZ^{1,n}/\spn{\omega_d}$, 
the $[\beta]$ are the weights in a single Weyl group orbit,
except for $d=7$ and $m=1,6$.
\end{enumerate}
These claims determine the required grading of the Lie algebra $\ee{8}$ from its root decomposition.

To verify (A), we need to know the roots of $\Etil_8$ in the range $1\leq m\leq d-1$, 
which are as follows (cf.~\cite[\S25.5.2]{Ma}), 
\begin{equation}\label{eq:enum-roots}
\begin{aligned}
 \text{(i)} &\quad\text{if $m=1$, then $\alpha = e_i-e_j$,} \\
 \text{(ii)} &\quad\text{if $1\leq m\leq 3$, then $\alpha=h-e_J$, with $|J|=3$,} \\ 
 \text{(iii)} &\quad\text{if $1\leq m\leq 6$, then $\alpha=2h-e_J$, with $|J|=6$,} \\ 
 \text{(iv)} &\quad\text{if  $m=d-1$,  then $\alpha = \omega_0 - e_i + e_j$,} \\ 
 &\text{where further $i\leq n < j$ and $J\not \sub \{1,\dots,n\}$.}
\end{aligned}
\end{equation}
Note that this way of enumerating the roots of $\Etil_8$ corresponds essentially to the decomposition 
$\ee{8}=\sl{9}\oplus \Lambda^3\oplus \Lambda^6$ from \eqref{eq:EMM}.

In each case, we see, as claimed, that $\gamma=-e_I$, where $|I|=m$.
If $n>1$, then any such $I$ can occur.
If $n=1$, then $m=4$ is not possible, and, if $n=0$, then only $m=3$ and $m=6$ are possible.

To verify (B), 
we first observe that $\ZZ^{1,n}/\spn{\omega_d}$ is also a known incarnation of the weight lattice of $\g_U$ 
and that the Weyl group action lifts to $\ZZ^{1,n}$, fixing $\omega_d$.
Then we can just use \eqref{eq:enum-roots} to explicitly enumerate the $\beta$ and
check they form a single orbit, except in the case that they don't
(see Appendix A for details).

We now repeat the analysis in the remaining case \eqref{eq:om_8b}, i.e.~$d=8b$,
by using the orthogonal direct sum decomposition
\[
   \ZZ^{1,9}= H \oplus \ZZ^8,
\]
where  $H=\langle f_1,f_2\rangle$ and $\ZZ^8=\langle e_3,\dots,e_9,e_0 \rangle$,
in terms of new elements 
\[
  f_1=h-e_1,\quad f_2=h-e_2, \quad e_0=h-(e_1+e_2).
\]
Note that $H$ is hyperbolic, i.e~$f_i^2=0$ and $f_1\cdot f_2=1$,
while $\ZZ^8$ is negative definite standard Euclidean.
As before, under this decomposition
$\omega_0=\omega_{8b} -\Delta_{8b}$, 
where $\omega_{8b}=2(f_1+f_2)$, from \eqref{eq:om_8b}, and $\Delta_{8b}=e_{3}+\dots+e_9+e_0$.

Also as before, write any root $\alpha=\beta+\gamma$,
set $m=\omega_{8b}\cdot\alpha=\omega_{8b}\cdot\beta=\Delta_{8b}\cdot\gamma$ 
and consider those $\Etil_8$ roots $\alpha$ with $0\leq m < 8$, noting that $m$ must be even.

If $m=0$ and $\alpha^2=-2$, then either it is a root of $\sl{8}$, 
whose root lattice is $\Delta_{8b}^\perp\sub\ZZ^8$,
or it is a root of $\g_U=\sl{2}$, 
whose root lattice is $\omega_{8b}^\perp\sub H$.

If $1\leq m\leq 7$, then explicit knowledge of the $\Etil_8$ roots gives
\begin{itemize}
 \item when $m=2$, we have $\alpha = f_1-e_I$ or $f_2-e_I$,
 \item when $m=4$, we have $\alpha=f_1+f_2-e_I$,
 \item when $m=6$,  we have $\alpha = 2f_1+f_2-e_I$ or $f_1+2f_2-e_I$
\end{itemize}
where, in each case, $I\sub\{3,\dots,9,0\}$ with $|I|=m$ and so $\gamma=-e_I$, as required.
Here we can see explicitly that $\beta^2=m-2$ and that the $[\beta]\in H/\langle \omega_{8b}\rangle$ 
are the weights of an irreducible representation of $\g_U=\sl{2}$, 
which is either fundamental, when $m=2,6$, or trivial, when $m=4$.

We can see the duality between $\gurep{m}{d}$ and $\gurep{d-m}{d}$
by observing that the map $\alpha\mapsto \omega_0-\alpha$ is an involution of
the set of $\Etil_8$ roots with $1\leq m\leq d-1$.
Under the decomposition $\alpha=\beta+\gamma$,
this becomes $\beta\mapsto\omega_d-\beta$ and $\gamma\mapsto\Delta_d-\gamma$,
so that $m\mapsto d-m$.
Projecting onto the respective weight lattices, these involutions are just
$[\beta]\mapsto-[\beta]$ and $[\gamma]\mapsto-[\gamma]$,
which always interchange the weights of a representation with those of its dual.
In our case, these representations are $\gurep{m}{d}$ and $\gurep{d-m}{d}$, 
or $\extpow{m}{d}$ and $\extpow{d-m}{d}$.
\end{proof}

For reference, Figure~\ref{fig:Rmd} records the dimensions of the representations $\gurep{m}{d}$,
summarising the calculations in Appendix A.

\begin{figure}[h]\centering
\newcommand{\xx}{0.9}
\begin{tikzpicture}[scale=0.75]
  \draw (-2.5*\xx,0) node {$d=2$};
  \draw (0,0) node {56};
  \draw (-3*\xx,-1) node {$d=3$};
   \foreach \m/\v in {1/27,2/27}
    {\draw (-1.5*\xx+\m*\xx,-1) node {\v};}
  \draw (-3.5*\xx,-2) node {$d=4$};
   \foreach \m/\v in {1/16,2/10,3/16}
    {\draw (-2*\xx+\m*\xx,-2) node {\v};}
  \draw (-4*\xx,-3) node {$d=5$};
    \foreach \m/\v in {1/10,2/5,3/5,4/10}
    {\draw (-2.5*\xx+\m*\xx,-3) node {\v};}
  \draw (-4.5*\xx,-4) node {$d=6$};
    \foreach \m/\v in {1/6,2/3,3/2,4/3,5/6}
    {\draw (-3*\xx+\m*\xx,-4) node {\v};}
  \draw (-5*\xx,-5) node {$d=7$};
    \foreach \m/\v in {1/3,2/2,3/1,4/1,5/2,6/3}
    {\draw (-3.5*\xx+\m*\xx,-5) node {\v};}
  \draw (-5.5*\xx,-6) node {$d=8a$};
    \foreach \m/\v in {1/1,2/1,3/1,4/0,5/1,6/1,7/1}
    {\draw (-4*\xx+\m*\xx,-6) node {\v};}
  \draw (-5.5*\xx,-7) node {$d=8b$};
    \foreach \m/\v in {1/0,2/2,3/0,4/1,5/0,6/2,7/0}
    {\draw (-4*\xx+\m*\xx,-7) node {\v};}
  \draw (-6*\xx,-8) node {$d=9$};
  \foreach \m/\v in {1/0,2/0,3/1,4/0,5/0,6/1,7/0,8/0}
    {\draw (-4.5*\xx+\m*\xx,-8) node {\v};}   
\end{tikzpicture}
\caption{Dimensions of $\gurep{m}{d}$, for $m=1,\ldots,d-1$.}
\label{fig:Rmd}
\end{figure}

\begin{remark} \label{rem:sugra}
Note that the decomposition in \eqref{eq:E8decomp} holds for both the complex form of $\ee{8}$
and the split real form.
The restriction of the Killing form gives a pairing between the degree $m$ and $d-m$ parts
of the grading, which must (up to scale) come from the dualities encountered in the proof.

To complete the supergravity picture, one would want to know that 
a split real Cartan involution for $\ee{8}$ may be given by combining 
\begin{itemize}
\item a standard Cartan involution for $\sl{d}$ coming from a metric, together with
 its Hodge-$*$ isomorphisms $\extpow{m}{d}\isom\extpow{d-m}{d}$,
\item  a split real Cartan involution of $\g_U$ and
 isomorphisms $\gurep{m}{d} \isom \gurep{d-m}{d}$ 
 compatible with the action of the compact part of $\g_U$. 
\end{itemize}

The lattices $\ZZ^{1,n}$ and $H$ that appear in the proof of Thm~\ref{thm:main}
can be identified as brane charge lattices as in \cite[\S3]{IqNeVa01}
and one should presumably expect that the lifts $\beta$ of the weights of $\gurep{m}{d}$
are precisely the charges of the $\half$-BPS branes in the relevant theory,
since these branes are charged by the form-fields appearing in supergravity multiplet
$(128)_B$.
\end{remark}

\begin{remark} \label{rem:dPlink}
The link with del Pezzo surfaces is also becoming apparent.
The lattice $\ZZ^{1,9}$ 
is the intersection lattice of a rational elliptic surface $X_0$ 
obtained by blowing up $\PP^2$ at 9 points 
and $\omega_0$ is its anti-canonical class.
The classes $\omega_d$ are the pullbacks of the anti-canonical classes of 
certain del Pezzo surfaces $X_d$ to which $X_0$ maps 
by blowing down the exceptional curves in classes $e_{n+1},\dots,e_9$.
The fundamental weights of $\Lambda(\Etil_8)$ are given by $\omega_1,\dots,\omega_6$,
as in \eqref{eq:om_d},
together with $\omega_{2'}=h-e_1$, $\omega_{3'}=h$ and $\omega_{4'}=2h-(e_1+e_2)$,
for the vertex labelling as in Fig.~\ref{fig:Etil-grading}.

The lattice $\ZZ^{1,n}$ is then 
naturally identified with the intersection lattice $I_X$ of a del Pezzo surface $X=X_{d}$ 
of degree~$d=9-n$, 
except in the case $d=8b$, when $I_X=H$.
In all cases, the weight lattice of $\g_U$ may be identified with the quotient $I_X/\spn{-K}$,
where $-K=\omega_d$ is the anti-canonical class.
\end{remark}

In fact, we can characterise the lifted weights of $\gurep{m}{d}$ as follows.

\begin{proposition}
\label{prop:beta-eqns}
The $\g_U$-weights of the representations $\gurep{m}{d}$ in Theorem~\ref{thm:main}
lift to the intersection lattice $I_X$ as the solutions $\beta$ of the equations
\begin{equation}
\label{eq:beta}
(i)\quad (-K)\cdot\beta = m,
\qquad
(ii)\quad \beta^2 =m-2.
\end{equation}
\end{proposition}

\begin{proof}
Note first that on $X_d$ the anti-canonical class $-K$ is $\omega_d$, so (i) is just the defining equation
$m=\omega_d\cdot\alpha=\omega_d\cdot\beta$, when $\beta$ is the $I_X$ component of an 
$\Etil_8$ root $\alpha$ with $0\leq m < d$.
Further, such $\alpha=\beta+\gamma$, where $\gamma$ is the $\ZZ^d$ component and 
we saw that $\gamma=-e_I=-\sum_{i\in I} e_i$, for $I$ an $m$-element subset of
$\{n+1,\dots,9\}$.
In particular, this means that $\gamma^2=-m$ and so (ii) is equivalent to $\alpha^2=-2$,
which holds precisely when $\alpha$ is a real root.
\end{proof}

Given \eqref{eq:beta}(i), the condition \eqref{eq:beta}(ii) is equivalent to $\beta\cdot(-K-\beta)=2$,
showing that $\beta\mapsto -K-\beta$ is an involution on the set of all lifted weights (i.e.~for all $m$) 
reflecting the duality $\gurepdual{m}{d}\isom \gurep{d-m}{d}$.

\newpage
\section{Rational curves and helical line bundles}\label{sec:RatExc}

We now explain how, on a Del Pezzo surface, the equations \eqref{eq:beta}
have a richer interpretation.

In general, the Riemann-Roch formula for a line bundle $L$ states that
\[
  \chi(L)= h^0(L)-h^1(L)+h^2(L) = \chi(\cO)+ \half L\cdot(L -K)
\]
and on a rational surface we have $h^1(\cO)=h^2(\cO)=0$, so $\chi(\cO)=1$.
Thus, if $L$ is a line bundle with $c_1(L)=\beta$, then the equations \eqref{eq:beta} 
from Proposition~\ref{prop:beta-eqns}
are equivalent to the equations
\begin{equation}
\label{eq:chi}
(i)\quad \chi(-L)=0,
 \qquad
(ii)\quad \chi(L)=m.
\end{equation}
On a del Pezzo surface, this condition has stronger consequences.

\begin{theorem}
\label{thm:excep}
Let $X$ be a del Pezzo surface of degree $d$ and $L$ be a line bundle on $X$,
other than $\cO$ or $-K$, satisfying
\eqref{eq:chi} with $0\leq m\leq d$. Then
\begin{equation}
\label{eq:cohvan}
 (i)\quad H^k(-L)=0, 
 \;\text{for all $k$,}
  \qquad 
 (ii)\quad
  H^k(L)=0,
  \;\text{for $k>0$.}
\end{equation}
In other words, the pair $(\cO,L)$ is \emph{strongly exceptional}.
Furthermore, so is the pair $(L,-K)$.
Finally $h^0(L)=m$ and, for $m>0$, 
a general section of $L$ vanishes on a smooth connected rational curve.
\end{theorem}

\begin{proof}
 Because $-K$ is ample, we know that, if $(-K)\cdot L<0$, then $H^0(L)=0$ 
 and indeed this is still true when $(-K)\cdot L=0$, provided $L\not\cong\cO$.
 Furthermore, by Serre duality, we have 
 \[
   H^2(L)\isom H^0(K -L)\dual.
 \] 
Note that $(-K)\cdot (K-L)=-d-m$ and $(-K)\cdot (K+L)=-d+m$, so in the range
$0\leq m\leq d$ and provided $L\not\cong\cO$ or $-K$, we deduce that
 \[
   H^2(L)=H^2(-L)=H^0(-L)=0.
\]
If $m=0$, then we also get $H^0(L)=0$ and hence also $H^1(L)=H^1(-L)=0$, using \eqref{eq:chi}.

If $m>0$, then we may still use \eqref{eq:chi} to deduce that $H^1(-L)=0$
and also that $h^0(L)>0$.
By Bertini's Theorem, a general section $s\in H^0(L)$ vanishes on a smooth curve $C$,
with $L=\cO(C)$.
Consider the two short exact sequences associated to multiplication by $s$,
\begin{gather*}
 0\to \cO(-C) \lra \cO \lra \cO_C \to 0 
 \\
  0\to \cO \lra \cO(C) \lra \cO_C(C) \to 0 
\end{gather*}
and their associated long exact sequences
\begin{gather*}
   \cdots\to H^{k}(\cO(-C)) \lra H^k(\cO) \lra H^k(\cO_C) \lra H^{k+1}(\cO(-C)) \to\cdots
 \\
   \cdots\to H^k(\cO) \lra H^k(\cO(C)) \lra H^{k}(\cO_C(C)) \lra H^{k+1}(\cO) \to\cdots
\end{gather*}
From the first, the vanishing of $H^{k}(\cO(-C))$, for all $k$, yields that $C$ is connected, i.e~$H^0(\cO)\to H^0(\cO_C)$ is an isomorphism, and that its genus $g_C=h^1(\cO_C)=0$, i.e~$C$ is rational.

From the second, we get $H^1(\cO(C)) \cong H^{1}(\cO_C(C))$.
Since the self-intersection $C^2=m-2\geq-1$, this implies that $H^1(L)=0$,
which completes the required cohomology vanishing for $L$.

Finally, by comparing with \eqref{eq:beta}, we note that $L'=-K-L$ satisfies the same conditions 
\eqref{eq:chi} as $L$ does,
but with $m'=d-m$, and so the same vanishing
holds for $L'$, that is,  $(L,-K)$ is also strongly exceptional. 
\end{proof}

\begin{remark}\label{rem:myst-dual}
The cohomology long exact sequences occurring in the proof of Theorem~\ref{thm:excep} can be used to show conversely that, for any smooth rational curve $C$, the line bundle $L=\cO(C)$ satisfies \eqref{eq:cohvan}. 
Thus the classes of smooth rational curves in the range 
$0< (-K)\cdot C < d$ are in one-one correspondence with the lifted weights
of the representations $\gurep{m}{d}$, and thus with $\half$-BPS brane charges,
as in the original formulation of mysterous duality \cite{IqNeVa01}.
The involution $L\leftrightarrow -K-L$ realises the duality  
$\gurepdual{m}{d} \isom \gurep{d-m}{d}$.
\end{remark}

\begin{remark}
If both $(\cO,L)$ and $(L,-K)$ are strongly exceptional,
then we say that $L$ is \emph{helical}.
This is a necessary (but not sufficient) condition for $L$ to appear in a geometric helix (\cite[Def. 1.5]{BS09}).
Informally, a helix provides a collection of `nice' bases for the derived category of $X$, 
or of the local Calabi-Yau $Y$ that is the total space of the canonical bundle $K_X$. 
More precisely, one seeks the period of a helix
\[
  \cO,L_1,\ldots,L_{n-1},-K
\]
in which all pairs except $(\cO,-K)$ are strongly exceptional.
The full helix of period length $n$ is the sequence $(L_k )_{k\in \ZZ}$ satisfying
$L_{k+n}=L_k - K$.
The first example of a full helix is the sequence $(\cO(k))_{k\in \ZZ}$ on $X=\PP^2$.
This has period $\cO,\cO(1),\cO(2),\cO(3)$ of length 3, because $-K=\cO(3)$.
\end{remark}

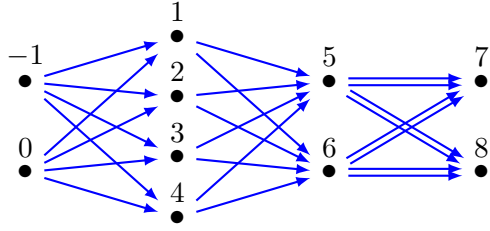
\begin{figure}[h]\centering
\begin{tikzpicture}[scale=0.8,
quiverarrow/.style={-latex, blue, thick}]
\foreach \n/\x/\y in {-1/-2.5/2.25, 0/-2.5/0.75, 1/0/3, 2/0/2, 3/0/1, 4/0/0, 5/2.5/2.25, 6/2.5/0.75, 7/5/2.25, 8/5/0.75}
 {\draw (\x,\y) node (v\n) [label={[label distance=-5pt] 90:{\small $\n$}}] {$\bullet$};}
\foreach \t/\h in {-1/1, -1/2, -1/3, -1/4, 0/1, 0/2, 0/3, 0/4, 1/5, 2/5, 3/5, 4/5, 1/6, 2/6, 3/6, 4/6}
 {\draw[quiverarrow] (v\t) -- (v\h);}
\foreach \t/\h in {5/7, 6/8}
 {\draw [transform canvas={yshift=0.25ex},quiverarrow] (v\t)--(v\h);
   \draw [transform canvas={yshift=-0.25ex},quiverarrow] (v\t)--(v\h);}
\draw [transform canvas={xshift=0.14ex, yshift=0.2ex}, quiverarrow] (v5)--(v8);
\draw [transform canvas={xshift=-0.14ex, yshift=-0.2ex}, quiverarrow] (v5)--(v8);
\draw [transform canvas={xshift=0.14ex, yshift=-0.2ex}, quiverarrow] (v6)--(v7);
\draw [transform canvas={xshift=-0.14ex, yshift=0.2ex}, quiverarrow] (v6)--(v7);
\end{tikzpicture}
\caption{A helix of period 8 on $dP_4$.}
\label{fig:helix}
\end{figure}

\begin{example}
A more interesting example is a helix of period length 8 on $dP_4$, depicted in Figure~\ref{fig:helix}.
A period is given by the line bundles
\begin{gather*}
L_0=\cO,\; L_1=\cO(e_0),\; L_2=\cO(e_3),\;  L_3=\cO(e_4),\;  L_4=\cO(e_5), \\
L_5=\cO(f_1),\; L_6=\cO(f_2),\; L_7=\cO(f_1+f_2),\;  L_8=-K
\end{gather*}
The 16 arrows in the first two columns correspond to the 16 exceptional lines, while the 4 double arrows in the third column correspond to 4 of the 10 pencils of conics on $dP_4$.
This choice of helix amounts to the choice of the four non-intersecting exceptional lines on $dP_4$ 
which are blown down to get $dP_{8b}=\PP^1\times\PP^1$.

At the level of Weyl groups and root systems, 
the choice of helix breaks the $\so{10}$ symmetry to 
$\so{4}\oplus\so{6}=\su{2}\oplus\su{2}\oplus\su{4}$.
This is the symmetry of the Pati-Salam model of fundamental particles \cite{PS74},
in which the 16 weights corresponding to the exceptional lines are the quantum numbers of the fermions.
The smaller Weyl group acts by permuting the vertices in a single column 
and the smaller root system is given by the differences between these vertices.
\end{example}

\clearpage
\section*{Appendix A}

We enumerate the weights of the $\g_U$ representations $\gurep{m}{d}$ using \eqref{eq:enum-roots}:
\[
\begin{tabular}{ccc|cccc|c}
 d & n & m & n & $\binom{n}{3-m}$ & $\binom{n}{6-m}$ & n & Total \\ \hline
 2 & 7 & 1 & 7 & 21 & 21 & 7 & 56 \\ \hline
 3 & 6 & 1 & 6 & 15 & 6 & & 27 \\
       & & 2 & & 6 & 15 & 6 & 27 \\ \hline
 4 & 5 & 1 & 5 & 10 & 1 & & 16\\
       & & 2 & & 5 & 5 & & 10 \\
       & & 3 & & 1 & 10 & 5 & 16\\ \hline
 5 & 4 & 1 & 4 & 6 & & & 10\\
       & & 2 & & 4 & 1 & & 5\\
       & & 3 & & 1 & 4 & & 5\\ 
       & & 4 & & & 6 & 4 & 10\\ \hline
 6 & 3 & 1 & 3 & 3 & & & 6\\
       & & 2 & & 3 & & & 3\\
       & & 3 & & 1 & 1 & & 2\\ 
       & & 4 & & & 3 & & 3\\ 
       & & 5 & & & 3 & 3 & 6\\ \hline
 7 & 2 & 1 & 2 & 1 & & & 3\\
       & & 2 & & 2 & & & 2\\
       & & 3 & & 1 & & & 1\\ 
       & & 4 & & & 1 & & 1\\ 
       & & 5 & & & 2 & & 2\\ 
       & & 6 & & & 1 & 2 & 3\\ \hline
 8 & 1 & 1 & 1 &&&& 1 \\
       & & 2 & & 1 &&& 1\\
       & & 3 & & 1 &&& 1\\ 
       & & 4 & &&&& 0\\ 
       & & 5 & && 1 && 1\\ 
       & & 6 & && 1 && 1\\ 
       & & 7 & &&& 1 & 1\\ \hline
\end{tabular}
\]
Each entry in the table counts a single orbit in $\ZZ^{1,n}$  for the subgroup $S_n$ of the Weyl group $W(\g_U)$,
consisting respectively of the elements of the form $e_i$, $h-e_J$, $2h-e_J$ or $\omega_d-e_i$,
for $i\in\{1,\ldots,n\}$ or $J\sub\{1,\ldots,n\}$ of the appropriate size.

For $n\geq 3$, there is an extra generator of $W(\g_U)$ given by the root reflection for $h-e_{123}$.
Then each row combines into a single orbit for $W(\g_U)$.
Hence, for $d\leq 6$, the representation $\gurep{m}{d}$ is irreducible and miniscule.
For $d\geq 7$ almost all rows have only one entry,
so $\gurep{m}{d}$ is still irreducible and miniscule, except for $\gurep{1}{7}$ and $\gurep{6}{7}$, 
which have two miniscule summands.

\clearpage

\end{document}